\documentclass[11pt,reqno]{amsart}
\usepackage{amssymb}
\usepackage{amsmath}
\usepackage{geometry}                
\usepackage{graphicx}
\usepackage{epstopdf}
\usepackage{mathrsfs}
\usepackage{color}
\usepackage{comment}
\newcommand{\Addresses}{{
  \bigskip
  \footnotesize

 \textsc{Department of Mathematics, KTH Royal Institute of Technology,
  Stockholm, Sweden}\par\nopagebreak
  \textit{E-mail addresses:} \quad \texttt{persj@kth.se, jostromb@kth.se}
}}

\title{Analysis of Schr\"odinger means}
\author{Per Sj\"olin and Jan-Olov Str\"omberg}

\renewcommand{\S}[0]{{\bf S}}

\newcommand{\R}[0]{{\mathbb R}}

\newcommand{\supp}[0]{\mbox{supp }}
\newcommand{\e}[0]{\mbox{e}}

\newtheorem{theorem}{Theorem}
\newtheorem{lemma}{Lemma}
\newtheorem{corollary}{Corollary}
\renewcommand{\smallskip}{}

\newcommand{\aname}{\em}
\newcommand{\jname}{}

\begin{document}
\begin{abstract}
We study integral estimates of maximal functions  for Schr\"odinger means.
\end{abstract}
\let\thefootnote\relax 
\maketitle
\section{Introduction}
For  $f\in L^2(\R^n), n\ge1$ and $a>0$ we set \[
\hat f(\xi)=\int_{\R^n} \e^{-i\xi\cdot x} f(x)\,dx, \xi \in \R^n,
\]
and \[ 
S_tf(x)=(2\pi)^{-n} \int_{\R^n} \e^{i\xi\cdot x} \e^{it|\xi|^a}\hat f(\xi)\, d\xi,\quad x\in\R^n, t\ge0.
\]
 For $a=2$ and $f$ belonging to the Schwartz class
$\mathscr{S}(\R^n)$ we set $u(x,t)=S_tf(x)$. It then follows that  $u(x,0)=f(x)$ and $u$ satisfies the Schr\"odinger equation 
$i\partial u/ \partial t=\Delta u$.
\par
We introduce Sobolev spaces $H_s=H_s(\R^n)$ by setting 
\[
H_s=\{f\in \mathscr{S}^\prime; \|f\|_{H_s} <\infty \}, s\in\R,
\]
where
\[
 \|f\|_{H_s}=\left(\int_{\R^n} (1+[\xi[^2)^s |\hat f(\xi)|^2\,d\xi \right)^{1/2}.
\]
Let $E$ denote a bounded set in $\R$. For $r>0$ we let $N_E(r)$ denote the minimal number $N$ of intervals $I_l, l=1,2,\dots N$, of length $r$
such that $E\subset \bigcup_1^NI_l$.
For $f\in\mathscr{S} $ we introduce the maximal function \[
S_E^*f(x)=\sup_{t\in E}\,|S_t f(x)|\,,\quad x\in\R^n.
  \]
  In Sj\"olin and Str\"omberg \cite{Sjo-Str}  we proved the following theorem.\\[.5cm]
{\bf Theorem A}  {\em Assume that $n\ge1$ and $a>0$ and $s>0$. if $f\in\mathscr{S} $ then one has \[
\int_{\R^n}|S_E^*f(x)|^2\,dx\lesssim \left(\sum_{m=0}^\infty N_E( 2^{-m})2^{-2ms/a} \right) \|f\|^2_{H_s}.
\]
}
\par
\vspace{.2cm}
Here we write $A\lesssim B$ if there is a constant $C$
such that $A\le C B$.\\
In the case $E=[0,1]$ it is easy to see that Theorem A implies the estimate
\begin{align}
\|S_E^*f\|_2\lesssim \|f\|_{H_s}
\end{align}
if $s>a/2$.\\
Let $(t_k)_1^\infty$ be a sequence satisfying
\begin{align}
 1>t_1>t_2>t_3>\dots>0\mbox{ and }\lim_{k\rightarrow\infty}t_k=0.
 \end{align}  
  Set \(
 A_j=\{ t_k; 2^{-j-1}<t_k\le 2^{-j} \}\mbox{ for }j=1,2,3,\dots .
 \)
Let $\#A$ denote the number of elements in a set $A$.\\
In \cite{Sjo-Str} we used Theorem A to obtain the following results.\\[.5cm]
{\bf Theorem B}  {\em Assume that $n\ge1$ and $a\ge2s$ and $0<s\le 1/2$, and $1/b>(a-2s)/ 2s$.\\Assume also that
\[
\#A_j\lesssim 2^{bj} \mbox{ for }j=1,2,3,\dots,
\]
and that $f\in H_s$. Then 
\begin{align}
\lim_{k\to\infty}S_{t_k} f(x)=f(x)
\end{align}
almost everywhere.\\[1cm]
}
{\bf Theorem C}  {\em Assume that $n\ge1$ and $a\ge2s$ and $0<s\le 1/2$, and that $\sum_1^\infty t_k^\gamma>\infty$, where
$1/\gamma>(a-2s)/2s$.\\
If also  $f\in H_s$  then (3) holds almost everywhere. \\
}
\par
\vspace{1cm}
Now let $E$ denote a bounded set in $\R^{n+1}$ and let $E_0$ denote the projection of $E$ onto the $t$-axis, i.e. 
$E_0=\{t;\mbox{ there extists }  x\in \R^n\mbox{ such that } (x,t)\in E\}$.\\
 For $a>0$ define an $a$-cube in $\R^{n+1}$
 with side $r>0$ as an axis-parallell rectangular box with sidelength  $r^a$ in the $t$-direction and with sidelength $r$ in the remaining directions
 $x_1,\dots,x_n$. Thus an $a$-cube in $\R^{n+1}$ has volume $r^{n+a}$.\\
 For $r>0$ let $N_{E,a}(r)$ denote the minimal number $N$ of $a$-cubes $Q_l, l=1,2,\dots,N$ of side $r$, such that $E\subset\bigcup_1^N Q_l $.\\
 \par
 For $f\in\mathscr{S} $ we introduce the maximal function \[
S_E^*f(x)=\sup_{(y,t)\in E}\,|S_t f(x+y)|\,,\quad x\in\R^n.
  \]
We shall prove the following inequality
\begin{theorem}
Assume that $n\ge1$ and $a>0$ and $s>0$. if $f\in\mathscr{S} $ then one has \[
\int_{\R^n}|S_E^*f(x)|^2\,dx\lesssim \left(\sum_{m=0}^\infty N_{E,a}( 2^{-m})2^{-2ms} \right) \|f\|^2_{H_s}.
\]
\end{theorem}
The following estimate follows directly.
\begin{corollary}
Assume $n\ge1, a>0, s>0$ and $f\in\mathscr{S} $ . If \[
\sum_{n=0}^\infty N_{E,a}( 2^{-m})2^{-2ms}<\infty,
\]
then
\begin{align}
\|S^*_Ef\|_2\lesssim\|f\|_{H_s}.
\end{align}
\end{corollary}
\par
Now let $\Gamma:[0,1]\to (\R^n)$  and let $E$ be a subset of the graph of  $\Gamma$ that is \[
E\subset\{(\Gamma(t),t); 0\le t\le1\}.
\]
 Assume that 
 \begin{align}
 |\Gamma(t_1)-\Gamma(t_2|\lesssim|t_1-t_2|^\beta\mbox{ for } t_1, t_2\in E_0,
 \end{align}
 where $0<\beta$.
Also set $a_1=1/\beta$ and $a_2=\max(a,a_1)$.\\
We shall first study the case when $E$ is the graph of $\Gamma$. Then (5) holds for any $\beta>1$ only if $\Gamma$ is constant.\\
It follows from a result of Cho, Lee, and Vargas \cite{Cho-Lee-Var} that if $n=1,a=2, B$ is an interval of $\R$ and $E$ is the graph of $\Gamma$, 
then \[
\|S_E^*f\|_{L^2(B)}\le C_B \|f\|_{H_s} \mbox{ if }s>\max(1/2-\beta,1/4), \mbox{ and }0<\beta\le1.
\]
We have the following result.
\begin{theorem}
Assume $n\ge1,a>0, f\in\mathscr{S} $ and that $E$ is the graph of $\Gamma$. Then (4) holds for $2s>a_2$. 
\end{theorem} 
\par 
We shall then study the case when $E=\{(\Gamma(t_k),t_k); k=1,2,3\dots\}$ where the sequence $(t_k)_1^\infty$ satisfies (2).  We  have the following results.
\begin{theorem}
Assume $n\ge1, a>0$, and  $f\in\mathscr{S} $.\\[.2cm]
In the case $2s>a_2$ then (4) holds.\\[.2cm]
In the case  $2s=a_2$ assume that $\sum_1^\infty t_k^\gamma<\infty$ for some $\gamma>0$. Then (4) holds.\\[.2cm]
In the case  $2s<a_2$ assume that $\sum_1^\infty t_k^\gamma<\infty$ for some $\gamma$ satisfying\\
 \(\gamma<2s/(a_2-2s).\) Then (4) holds.
 \end{theorem}
\section{Proof of Theorem 1}
If $y$ and $y_0$ belong to to $\R^n$ we write $y=(y_1,\dots,y_n)$ and $y_0=(y_{0,1},\dots,y_{0,n})$. We shall give the proof of Theorem 1.
\begin{proof}
Assume $y_0\in\R^n,t_0\in\R, 0<r\le1$ and let \[
E=\{(y,t)\in\R^{n+1}; y_{0,j}\le y_j\le y_{0,j}+r \mbox{ for } 1\le j\le n,\mbox{ and }t_0\le t\le  t_0 +r^a \}
\]
We have \[
S_tf(x+y)=c\int \e^{i\xi\cdot x} \e^{i\xi_1y_1}\dots\e^{i\xi_ny_n} \e^{it|\xi|^a} \hat f(\xi)\, d\xi
\]
where $c=(2\pi)^{-n}$, and for $1\le j\le n$ we write $\e^{i\xi_jy_j}=\Delta_j+ \e^{i\xi_jy_{0,j}}$  where \[
\Delta_j= \e^{i\xi_jy_j} -\e^{i\xi_jy_{0,j}}.
\]
We also write $\e^{it|\xi|^a}=\Delta_{n+1}+\e^{it_0|\xi|^a}$ where \[
\Delta_{n+1}=\e^{it|\xi|^a}-\e^{it_0|\xi|^a}.
\]
Hence \[
S_tf(x+y)=c\int\e^{i\xi\cdot x}\left(\Delta_1+\e^{i\xi_1y_{0,1}} \right)\dots\left(\Delta_n+\e^{i\xi_ny_{0,n}} \right)\\
\left(\Delta_{n+1}+\e^{it_0|\xi|^a} \right)\hat f(\xi)\, d\xi,
\]
and it follows that $S_tf(x+y)$ is the sum of integrals of the form 
\begin{align}
c\int\e^{i\xi\cdot x}\left(\prod_{j\in D} \Delta_j  \right) \left(\prod_{j\in B} \e^{i\xi_jy_{0,j}} \right)\Delta_{n+1} \hat f(\xi)\,d\xi,
\end{align}
or
\begin{align}
c\int\e^{i\xi\cdot x}\left(\prod_{j\in D} \Delta_j  \right) \left(\prod_{j\in B} \e^{i\xi_jy_{0,j}} \right)\e^{it_0|\xi|^a} \hat f(\xi)\,d\xi,
\end{align}
Here $D$ and $B$ ar disjoint subsets of $\{1,2,3,\dots,n\}$ and $D\cup B=\{1,2,3,\dots,n\}$. We denote the integrals in (6) 
 by $S^\prime_tf(x,y)$ and shall describe how they can be estimated. The same argument works also for integrals in (7).\\[.1cm]
\par
For $j\in D$ we write\[
\Delta_j=i\xi_j\int_{y_{0,j}}^{y_j} \e^{i\xi_js_j}\,ds_j,
\] 
and we also write\[
\Delta_{n+1}=i|\xi|^a\int_{t_0}^t\e^{i|\xi|^as_{n+1}}\,ds_{n+1}.
 \]
 Assuming $D=\{k_1,k_2,\dots,k_p\}$ we then have \[
 \begin{array}{ll}
 S^\prime_tf(x,y)=\int_{\R^n} \int_{y_{0,k_1}}^{y_{k_1}}\int_{y_{0,k_2}}^{y_{k_2}}\dots\int_{y_{0,k_p}}^{y_{k_p}}\int_{t_0}^t &\e^{i\xi\cdot x}
 \left(\prod_{j\in D} i\xi_j\e^{i\xi_js_j}\right) \left(\prod_{j\in B} \e^{i\xi_jy_{0,j}}\right)\\[.2cm]
 &\cdot i|\xi|^a\e^{i|\xi|^as_{n+1}}\hat f(\xi)\,ds_{k_1}ds_{k_2}\dots ds_{k_p}ds_{n+1}d\xi.
 \end{array}
 \]
 Changing the order of integration one then obtains \[
 \begin{array}{ll}
 |S^\prime_tf(x,y)|\le c\int_{y_{0,k_1}}^{y_{k_1}}\int_{y_{0,k_2}}^{y_{k_2}}\dots\int_{y_{0,k_p}}^{y_{k_p}}\int_{t_0}^t &
 \left| \int_{\R^n}\e^{i\xi\cdot x}
 \left(\prod_{j\in D} i\xi_j\e^{i\xi_js_j}\right) \left(\prod_{j\in B} \e^{i\xi_jy_{0,j}}\right)\right.\\[.2cm]
 &\left.\cdot i|\xi|^a\e^{i|\xi|^as_{n+1}}\hat f(\xi)\,d\xi\right|\,ds_{k_1}ds_{k_2}\dots ds_{k_p}ds_{n+1},
\end{array}
\]
or
 \[
 S^\prime_tf(x,y)=\int_{y_{0,k_1}}^{y_{k_1}}\dots\int_{y_{0,k_p}}^{y_{k_p}}\int_{t_0}^t F_D(x;s_{k_1},,\dots ,s_{k_p},s_{n+1})\,ds_{k_1}ds_{k_2}\dots ds_{k_p}ds_{n+1},
 \]
 where\[
 F_D(x;s_{k_1},s_{k_2},\dots,s_{k_p},s_{n+1})=
c\int_{\R^n}\e^{i\xi\cdot x}
 \left(\prod_{j\in D} i\xi_j\e^{i\xi_js_j}\right) \left(\prod_{j\in B} \e^{i\xi_jy_{0,j}}\right)
i|\xi|^a\e^{i|\xi|^as_{n+1}}\hat f(\xi)\,d\xi.
\] 
and\[
\begin{array}{l}
\hspace{-.3cm}\sup\limits_{\quad(y,t)\in E}|S^\prime_tf(x,y)|
\\[.2cm]\hspace{.5cm}
\le\int_{y_{0,k_1}}^{y_{0,k_1}+r}\dots\int_{y_{0,k_p}}^{y_{0,k_p}+r}\int_{t_0}^{t_0+r^a} 
|F_D(x;s_{k_1},\dots ,s_{k_p},s_{n+1})|\,ds_{k_1}\dots ds_{k_p}ds_{n+1}.
\end{array}
\]
Invoking Minkowski's inequality and Plancherel's formula we then obtain\[
\begin{array}{l}
\left(\int_{\R^n}\hspace{-.3cm}\sup\limits_{\quad(y,t)\in E}|S^\prime_tf(x,y)|^2\,dx\right)^{1/2}\\[.5cm]
\hspace{.3cm}\le\int_{y_{0,k_1}}^{y_{0,k_1}+r}\dots\int_{y_{0,k_p}}^{y_{0,k_p}+r}\int_{t_0}^{t_0+r^a}
 \|F_D(\cdot;s_{k_1},\dots ,s_{k_p},s_{n+1})\|_2\,ds_{k_1}\dots ds_{k_p}ds_{n+1}\\[.5cm]
 \hspace{.3cm}=c(2\pi)^{n/2}\int_{y_{0,k_1}}^{y_{0,k_1}+r}\dots\int_{y_{0,k_p}}^{y_{0,k_p}+r}\int_{t_0}^{t_0+r^a}
 \left(\int_{\R^n}(\prod\limits_{j\in D}|\xi_j|^2)|\xi|^{2a}|\hat f(\xi)|^2d\xi\right)^{1/2}ds_{k_1}\dots ds_{k_p}ds_{n+1}\\[.5cm]
  \hspace{.3cm}\le r^pr^aA^pA^a \|f\|_2 
\end{array}
\]
if $f\in L^2(\R^n))$ and $\supp\hat f\subset B(0,A)$, where $A\ge1$.\\
With similar arguments we get the estimate  $r^pA^p \|f\|_2$ for the integrals in (7), and by summation of the all integrals of forms
 (6) or (7) we get
 \begin{align}
 \|\hspace{-.3cm}\sup\limits_{\quad(y,t)\in E}|S_tf(x+y)|\|_2\le (1+rA)^n(1+r^aA^a)\|f\|_2
 \end{align}
 \par
 Now let $E$ be a set in $\R^{n+1}$ with the property that $E\subset\bigcup_1^NQ_l$ where the sets $Q_l$ are $a$-cubes
 with side $r$ with $rA\le1$ of the type we have just considered. One has \[
 \sup\limits_{\quad(y,t)\in E}|\S_tf(x+y)|^2\le\sum\limits_1^N\sup\limits_{(y,t)\in Q_l}|S_tf(x+y)|^2
 \]
 and it follows from (8) that
 \begin{align}
 \int_{\R^n}\hspace{-.3cm}\sup\limits_{\quad(y,t)\in E}|S_tf(x+y)|^2\,dx\le 2^{2n+2}N\|f\|_2^2
 \end{align}
if $rA\le1$ and $\supp\hat f\subset B(0,A)$.\\
Now let  $f\in\mathscr{S}$. We write $f=\sum_{k=0}^\infty f_k$, where the functions $f_k$ are defined in the following way. We set
$\hat f_0(\xi)=\hat f(\xi)$ for $|\xi|\le1$  and $\hat f_0(\xi)=0$ for $|\xi|>1$. For $k\ge1$ we let $\hat f_k(\xi)=\hat f(\xi)$ for
$2^{k-1}<|\xi|\le2^k$ and $\hat f_k(\xi)=0$ otherwise.\\
Choosing  real numbers $g_k>0, k=0,1,2,\dots,$ we have\[
\begin{array}{l}
\sup\limits_{\quad(y,t)\in E} |S_t f(x+y)|\le\sum\limits_{k=0}^\infty \sup\limits_{\quad(y,t)\in E} |S_t f_k(x+y)|\\[.5cm]
\quad=\sum\limits_{k=0}^\infty  g_k^{-1/2}\sup\limits_{\quad(y,t)\in E} |S_t f_k(x+y)|\,g_k^{1/2}\\[.5cm]
\le \left(\sum\limits_{k=0}^\infty  g_k^{-1}\sup\limits_{\quad(y,t)\in E} |S_t f_k(x+y)|^2 \right)^{1/2}  \left(\sum\limits_{k=0}^\infty  g_k\right)^{1/2}
\end{array}
\] 
and invoking inquality (9) we also have \[
\begin{array}{l}
\int\hspace{-.3cm}\sup\limits_{\quad(y,t)\in E} |S_t f(x+y)|^2\,dx \\[.5cm]
\hspace{.5cm}\le \left(\sum\limits_{k=0}^\infty  g_k^{-1}\int\hspace{-.3cm}\sup\limits_{\quad(y,t)\in E} |S_t f_k(x+y))|^2 \,dx\right)  \left(\sum\limits_{k=0}^\infty  g_k\right)
\le\left(\sum\limits_{k=0}^\infty  g_k\right)\left(\sum\limits_{k=0}^\infty  g_k^{-1}\, 2^{2n+2}N_{E,a}(2^{-k})\|f_k\|_2^2\right).
\end{array}
\]
Chosing $g_k=N_{E,a}(2^{-k})2^{-2ks}$ we conclude that \[
\begin{array}{ll}
\int\hspace{-.3cm}\sup\limits_{\quad(y,t)\in E} |S_t f(x+y)|^2\,dx
&\le2^{2n+2}\left( \sum\limits_{k=0}^\infty N_{E,a}(2^{-k})2^{-2ks}\right) \left( \sum\limits_{k=0}^\infty2^{2ks}\|f_k\|_2\right)\\[.5cm]
&\lesssim\left( \sum\limits_{k=0}^\infty N_{E,a}(2^{-k})2^{-2ks}\right)\|f\|_{H_s}^2,
\end{array}
\]
and the proof of Theorem 1 is complete.
\end{proof}
\section{Proof of Theorems 2 and 3}
The following two lemmas follow easily from the definition of $N_{E,a}(r)$.
\begin{lemma}
Assume tha $0<r\le1$ and $0<b<b_1$. Then\[
N_{E,b}(r)  \le N_{E,b_1}(r),
\]
and\[
N_{E,b_1}(r) \le r^{b-b_1} N_{E,b}(r).
\]
\end{lemma}
\begin{lemma}
Assume that $E$ is a subset of the graph of $\Gamma$ satisfying (5)  with $\beta>0$
and let $b\ge1/\beta$\\
Then \[
N_{E,b}(r)\lesssim  N_{E_0}(r^b) \mbox{ for } 0<r\le1. 
\]
 ($E_0$, $\Gamma,\beta$ are defined  and equation (5) is in Section 1).
\end{lemma}
We shall then give the proofs of Theorems 2 and 3.
\begin{proof}[Proof of Theorem 2]
We have $a\le a_2$. Invoking Lemma 1 and Lemma 2 we obtain\[
N_{E,a}(2^{-m})\le N_{E,a_2}(2^{-m})\lesssim N_{E_0}(2^{-ma_2})\le 2^{ma_2}+1
\]
and\[
\sum\limits_{k=0}^\infty N_{E,a}(2^{-m})2^{-2ms}<\infty
\]
if $a_2-2s<0$, i.e.$s>a_2/2$.\\
Using Corollary 1 we condclude that\[ 
\|S^*_Ef\|_2\lesssim \|f\|_{H_s}
\]
if $s>a_2/2$.
\end{proof}
We shall then prove Theorem 3.
\begin{proof}[Proof of Theorem 3]
In the case $2s>a_2$ we can use the same argument as in the proof of Theorem 2 to prove that (4) holds.\\
Then assume $2s\le a_2$.
We have $E_0=\{t_k;k=1,2,3,\dots\}$ and assuming $\sum_1^\infty t_k^\gamma<\infty$ one obtains $\#A_j\lesssim2^{\gamma j}$. It then follows from 
Lemma 6 in \cite{Sjo-Str} that $N_{E_0}(2^{-m})\lesssim2^{\gamma m/(\gamma+1)}$.\\
We  have  $a\le a_2$. Applying Lemma 1 and Lemma 2 one obtains\[
N_{E,a}(2^{-m})\le N_{E,a_2}(2^{-m})\lesssim N_{E_0}(2^{-ma_2})\lesssim2^{a_2\gamma m/(\gamma+1)}.
\]
It follows that\[
 \sum\limits_{k=0}^\infty N_{E,a}(2^{-m})2^{-2ms}<\infty
\] 
if $\gamma a_2/(\gamma+1)<2s$, that is $a_2\gamma<2s\gamma+2s$ . \\
In the case $2s=a_2$ this holds for every $\gamma>0$.\\
 In the case $2s<a_2$ one has to assume $\gamma< 2s/(a_2-2s)$.\\[.2cm]
This completes the proof of Theorem 3.
\end{proof}

\Addresses
\end{document}